\documentclass[12pt,oneside,a4paper]{article}
\usepackage{makeidx}     
\usepackage{enumerate} 
\usepackage{graphicx,epsfig,lscape,float}    
\usepackage{amsmath,array,amssymb, amsfonts, latexsym, mathrsfs,amsthm}
\usepackage{mathtools}
\makeatother

\usepackage{calc}
\usepackage{ifthen}
\usepackage{xspace}

\usepackage{datapie}
\usepackage{adjustbox}

\usepackage{tikz}
\usetikzlibrary{positioning}
\usepackage{varwidth}

\usepackage{psfrag,fancybox}
\usepackage{mathtools}
\usepackage{tikz}
\usepackage{verbatim}
\usetikzlibrary{shapes,decorations}
\usetikzlibrary{arrows}
\usepackage{rotating}
\usepackage{pgfplots}
\usepgfplotslibrary{units}

\makeindex             
\usepackage{csquotes}

\newcommand{\restr}[1]{\lower0.4ex\hbox{$\vert$}\lower0.7ex\hbox{ $\!_{#1}$ }}

\newcommand{\R}{\mathbb{R}}

\newcommand{\Ga}{\Gamma}
\newcommand{\Om}{\Omega}

\newcommand{\bG}{\bar{\Ga}}

\newcommand{\conv}{\operatorname{conv}}

\newcommand{\define}{\coloneqq}

\newtheorem{theorem}{Theorem}[section]

\newtheorem{observation}[theorem]{Observation} 

\usepackage{amsmath,array,amssymb,amsfonts,graphicx}

\pgfplotsset{compat=1.18}

\begin{document}

 \title{Edge-unfolding nested prismatoids}

\author{Manuel Radons \\ radons@math.tu-berlin.de}
\date{Technische Universität Berlin, Chair of Discrete Mathematics/Geometry\\
}

\maketitle

\textbf{Abstract:}
{\em 
	A $3$-prismatoid is the convex hull of two convex polygons $A$ and $B$ which lie
	in parallel planes $H_A, H_B\subset\mathbb{R}^3$. Let $\tilde{A}$ be the orthogonal projection 
	of $A$ onto $H_B$. A $3$-prismatoid is called  nested if $\tilde{A}$ is properly contained in $B$, or vice versa. 
	We show that every nested $3$-prismatoid has an edge-unfolding to a non-overlapping polygon in the plane.
   \em} 
\section{Introduction}

Throughout, a \emph{polytope} is the boundary of a convex $3$-dimensional polytope.
If a polytope is cut along a spanning tree of its edge-graph, the resulting ($2$-dimensional) polyhedral surface can be \emph{unfolded} flat into the plane, cf. \cite{Ghomi2014Affine}. 
This unfolded surface, or: \emph{unfolding}, is called simple if it has no self-overlaps.
If a polytope has a simple unfolding, it is called \emph{unfoldable}.
The question whether every polytope is unfoldable can be dated back to the \enquote{Painter's Manual} by Albrecht D\"urer \cite{Duerer1525manual}.
It is thus often referred to as D\"urer's Problem.

It was proved by Ghomi that every polytope is unfoldable after an affine stretching, which implies that every combinatorial type of polytope has an infinite number of unfoldable realizations \cite{Ghomi2014Affine}.
O'Rourke established the unfoldability of nearly flat, acutely triangulated \emph{convex caps} \cite{Orourke2018flatcaps,Orourke2017flatcaps_add}.
A convex cap is a polytope $C$ which has a designated facet $F$ so that the orthogonal projection of $C\setminus F$ onto $F$ is one-to-one. 
An acute triangulation of a $2$-dimensional surface is a triangulation of the latter so that every interior angle of every of its triangles is smaller than $\pi/2$.  
A recent negative result, which Barvinok and Ghomi distilled from a highly original but flawed preprint of Tarasov \cite{Barvinok2017PseudoEdgeUO,Tarasov2008ExistenceOA}, concerns the existence of counterexamples to a more general
form of D\"urer's problem which considers cuts along so-called pseudo-edges, which are geodesics in the intrinsic metric of
a polytope. 
Another generalized form of D\"urer's problem concerns unfoldability of non-convex polytopes which are
combinatorially equivalent to a convex $3$-polytope.
There are several ununfoldable families of such polytopes known, cf. \cite{gruenbaum2002, Tarasov, DDE2020}. 

\subsection{Unfolding prismatoids}
In the context of this work, a \emph{prismatoid} $P$ is boundary of the convex hull of two convex polygons $A$ and $B$ that lie in parallel planes, say $H_A$ and $H_B$. As such, prismatoids are $3$-polytopes. We will refer to $A$ and $B$ as the \emph{top} and \emph{base} of $P$, respectively.
 The set of lateral facets of a prismatoid is called its \emph{band}.
 A prismatoid whose lateral facets are trapezoids is called a \emph{prismoid}. 
 In this case, corresponding top and base edges are parallel.
  
  \begin{figure}
  	\begin{center}
  		\begin{tikzpicture}
  		
  		\draw[thick] (2,2) -- (2,0.6) -- (3,0.6) -- (3,2) -- cycle; 
  		\draw[thick] (2,2) -- (2.5,2.5) -- (3,2);
  		\draw[thick] (2,0.6) -- (2.5,0.1) -- (3,0.6);
  		\draw[thick] (2,2) -- (1,3) -- (1.5,3.5) -- (2.5,2.5);
  		\draw[thick] (3,2) -- (4,3) -- (3.5,3.5) -- (2.5,2.5);
  		
  		\draw[thick] (6,2) -- (8.5,2) -- (8.5,0.6) -- (6,0.6) -- cycle; 
  		\draw[thick] (6,0.6) -- (6.5,0.1) -- (7,0.6) -- (7,2) -- (6.5,2.5) -- (6,2);
  		\draw[thick] (7.75,0.6) -- (7.75,2);
  	
  		\end{tikzpicture}
  	\end{center}
  	\caption{Petal and band unfolding of a prism over a triangle.}\label{fig:patterns}
  \end{figure}
  
  There are two natural ways to unfold prismatoids, the \emph{band unfolding}, and the \emph{petal unfolding} \cite{Orourke2012patches}, cf. Figure \ref{fig:patterns}.
  In a band unfolding the top and base are removed, and one lateral edge is cut. 
  Then the band is unrolled into the plane as one connected patch, and $A$ and $B$ are reattached to the latter along one suitable edge each. 
  Every prismoid has a lateral edge $e$ so that the band, if cut along $e$, can be unfolded without self-intersections \cite{Aloupis2005phd, aloupis2008paperversion}.
  But there exist prismoids, and hence prismatoids, so that every placement of the prismoid top overlaps with every one of its band unfoldings \cite{Orourke2007band_prismatoids}.
  
  In a petal unfolding either $A$ or $B$ is a designated facet to which all lateral facets are left attached. Assume that the designated facet is $B$. 
  Then for each vertex $b_i$ of $B$ exactly one lateral edge adjacent to it is cut. 
  The so-called resulting \emph{petals} are unfolded into the plane while $A$ is left attached to this unit along a single suitable edge.
  O'Rourke proved that every prismoid has a non-overlapping petal unfolding \cite{Orourke2001prismoids}.
  Smooth prismatoids, which are the convex hull of two smooth convex curves lying in parallel planes, have a petal unfolding as well \cite{Orourke2004smooth_prismatoids}.
  Further, several subclasses of prismatoids are known to have a petal unfolding. 
  A nonobtuse triangle is a triangle so that all its interior angles are smaller than or equal to $\pi/2$.
  A prismatoid has a petal unfolding if all its facets, except possibly its base $B$, are nonobtuse triangles \cite{Orourke2012patches}, or if the base is a rectangle and all other facets are acute triangles, or if the top and base are sufficiently far from each other \cite{bian2021prismatoids}. 
  
  \subsection{Main result}
  
  A prismatoid is called \emph{nested} if the orthogonal projection of (the convex polygon) $A$ onto $H_B$ is properly contained in (the convex polygon) $B$, or vice versa.
  (Note that nested prismatoids are a special class of convex caps.)
  Our main result is the following:
  \begin{theorem}\label{thm:nested}
  	Let $P$ be a nested prismatoid.
  	Then $P$ is unfoldable.
  \end{theorem}
  
  To this end we apply a combination of the petal and the band unfolding strategies to nested prismatoids.
  More precisely, we cut the band into two pieces.
  Crucial in the selection of the band-patches which are left intact is the notion of radially monotone polygonal paths, which was introduced and exploited to great effect in \cite{Orourke2018flatcaps}.

  \subsection{Content and Structure}
In Section \ref{sec:prelim} we will introduce the necessary concepts for our investigation and some preliminary results. 
Sections \ref{sec:cutting} and \ref{sec:gluing} contain the proof of our main result. 


\section{Preliminaries}\label{sec:prelim}

We will follow notation and naming conventions established
by Ghomi in \cite{Ghomi2014Affine} and O'Rourke  in \cite{Orourke2018flatcaps, aloupis2008paperversion} whenever possible.
As noted above, by \enquote{polytope} we mean the boundary of a convex $3$-polytope which lives in $\R^3$.
A prismatoid $P$ is the boundary of the convex hull of two convex polygons $A\subset H_A$ and $B\subset H_B$, where $H_A$ and $H_B$ are parallel affine planes. 
We will assume without loss of generality that the projection of $A$, the top of $P$, to $H_B$ is properly contained in $B$, the base of $P$. 
Further, we will assume that $H_B$ is the $xy$-plane embedded in $\R^3$ and $H_A$, which is a parallel to $H_B$, has a positive height.
Vertices of $A$ and $B$ are denoted $a_i$ and $b_j$, respectively. 

Vertices of an $n$-gon are enumerated counterclockwise from $1$ to $n$ with respect to a viewpoint above the polygon.
Let $D$ be a convex polygon bounded by the polygonal path $[d_1,\dots, d_n]$.
A subpath $[d_i, d_{i+1}, \dots, d_j]$ of the boundary of $D$ is denoted $(d_i,d_j)$.
We define the \emph{curvature at a vertex} $d_k$ as the angle spanned by the outward edge normals of $D$ at $d_k$. 
The \emph{total curvature of a subpath} $(d_i,d_j)$ is the sum of the curvatures at its interior vertices $\{d_{i+1}, d_{i+2},\dots,d_{j-1}\}$.
In some sources these quantities are also called the \emph{turn angle} and the \emph{turn}; see, for example, \cite{Aloupis2005phd}.
For brevity we will refer to the total curvature of a path simply as its curvature, if the meaning is clear from the context. 

A \emph{band piece} of a prismatoid is a connected set of its lateral facets. If a band piece is bounded at the top by
$\Ga\define [a_k,\dots,a_\ell]$ and at the base by $\Omega\define [b_K,\dots,b_L]$, then we call $\Ga$ its top boundary and the curvature of $\Ga$ its top curvature. Base boundary and curvature are defined analogously.
We use the same letters in upper and lower case to underline the facts that the end vertices of top and base boundary a) correspond to each other, while b) they usually do not have the same index.

We define the \emph{flat prismatoid} $P^0$ corresponding to $P$ as follows: 
The lower facet of $P^0$ coincides with the base $B$ of $P$.
The upper facets of $P^0$ are obtained as the cells of a subdivision 
of $B$ which is induced by the orthogonal projection of $P\setminus B$ onto $B$. This projection is one-to-one since nested prismatoids are polyhedral caps.

\subsection{Projections and unfoldings}

Let $P$ be a nested prismatoid with top $A$ and base $B$.
Then we denote the orthogonal projection of any subset $C\subset P$ onto $H_B$ by $\tilde{C}$.

The edge graph of a polytope is the graph whose vertex set is the set of vertices of the polytope, and whose edge set is the set of edges of the polytope.
As mentioned in the introduction, any spanning tree of the edge graph of a polytope $P$ induces an unfolding of $P$ into the plane, cf. \cite{Ghomi2014Affine}.
This unfolding is an \emph{isometric immersion} which is unique up to rigid motions.
The immersion is one-to-one, i.e., an embedding, if and only if the unfolding is \emph{simple} in the sense that it does not self-overlap.
In this case, and this case only, the image of the edge graph spanning tree, along which the polytope has been cut, is a simple -- usually non-convex -- polygon.

Any subset $C$ of a polytope $P$ which has been cut along a spanning tree of its edge graph can be isometrically immersed, i.e., unfolded into a plane, as well. 
We denote the unfolded image of $C$ by $\bar{C}$ and assume the plane to be $H_B$, the $xy$-plane embedded in $\R^3$. 
In particular, $\bar{P}$ denotes the unfolded image of $P$. 
In some sources $\bar{P}$ is subscripted with a $T$, where $T$ denotes the specific edge graph spanning tree along which $P$ has been cut.
But we will construct only one such tree. Hence, omitting the subscript introduces no ambiguities.



\subsection{Radial monotonicity}
A polygonal path $\Ga\subset\R^2$ is called \emph{radially monotone} if when traversing it from an endpoint the Euclidean distance to that endpoint increases monotonically.
In other words, for points $p_i,p_j,p_k\in\Ga$ the Euclidean distance of $p_i$ and $p_j$ is greater than that of $p_i$ and $p_k$ if and only if their distances within $\Ga$ have the same relation; cf. \cite{Orourke2018flatcaps}.
Let $\Gamma\define(d_i,d_j)$ be a subpath of a convex $n$-gon $D\define\conv(d_1,\dots,d_n)$ in some plane $H$. Furhter, let $\Om$ be another polygonal path in $H$ that is obtained from $\Gamma$ by 
decreasing the curvature at all of its interior vertices, but not decreasing curvature to $0$ or less at any vertex.
We then say that $\Om$ is obtained by \emph{stretching} $\Ga$.
In \cite{Orourke2018flatcaps} we find the following crucial observation. 
\begin{observation}\label{obs:monotone-intersections}
Let $\Gamma\define(d_i,d_j)$ be a subpath of a convex $n$-gon $D$ in some plane $H$. Further, let $\Om$ be obtained by stretching $\Ga$, keeping its first edge, $(d_i,d_{i+1})$, fixed. 
If $(d_{i+1},d_j)$ is radially monotone, then the intersection of $\Gamma$ and $\Om$ consists of the fixed edge $(d_i,d_{i+1})$ and no other point of the two paths, cf. Figure \ref{fig:monotone}. 
\end{observation}
Note that for a subpath of the boundary of a convex polygon a sufficient condition for radial monotonicity is that its curvature does not exceed $\pi/2$.
Proof: Assume there are points $p_1, p_2$ with the same distance to the starting point $p$ of a path $\Ga$. Then the equilateral triangle $[p,p_1,p_2,p]$ has interior angles $\alpha_1=\alpha_2<\pi/2$ at $p_1$ and $p_2$, respectively.
But by construction $\pi-\alpha_1=\pi-\alpha_2$---the curvature of the triangle triangle $[p,p_1,p_2,p]$ at $p_1$ and $p_2$, respectively---is a lower bound for the curvature of $\Ga$.

\begin{figure}
	\begin{center}
		\begin{tikzpicture}[scale=.9]
		\draw[thick,color=lightgray] (0,0) -- (-1,0);
		\draw[thick,color=lightgray] (0,0) -- (3,.5);
		\draw[thick,color=lightgray] (3,.5) -- (5,1.5);
		\draw[thick,color=lightgray] (5,1.5) -- (3,2.5);
		
		\draw[thick] (0,0) -- (-1,0);
		\draw[thick] (0,0) -- (2.9,.7);
		\draw[thick] (2.9,.7) -- (4.9,1.9);
		\draw[thick] (4.9,1.9) -- (2.7,2.5);
		\draw (3,1.3) node{$\Ga$};
		\draw[color=lightgray] (4,0.5) node{$\Om$};
		\draw (-1,-0.4) node{$d_i$};
		\draw (0.2,-0.41) node{$d_{i+1}$};
		\end{tikzpicture}
		\begin{tikzpicture}[scale=.9]
		\draw[thick,color=lightgray] (0,0) -- (3,.5);
		\draw[thick,color=lightgray] (3,.5) -- (5,1.5);
		\draw[thick,color=lightgray] (5,1.5) -- (5.7,2);
		\draw (3,1.3) node{$\Ga$};
		\draw[color=lightgray] (4,0.5) node{$\Om$};
		\draw (-1,-0.4) node{$d_i$};
		\draw (0.2,-0.41) node{$d_{i+1}$};
		\draw[thick] (0,0) -- (-1,0);
		\draw[thick] (0,0) -- (2.9,.7);
		\draw[thick] (2.9,.7) -- (4.9,1.9);
		\draw[thick] (4.9,1.9) -- (5.5,2.5);
		\draw (6,0) node{};
		\end{tikzpicture}
	\end{center}
	\caption{Stretching of paths with curvature $>\frac\pi2$ (left), resp. $\leq\frac\pi2$ (right).}\label{fig:monotone}
\end{figure}

\subsection{Flattening of the band}
The observation below was stated for the special case of nested prismoids in both \cite{Orourke2007band_prismatoids} and \cite{aloupis2008paperversion}.
We will generalize it to nested prismatoids.
\begin{observation}\label{obs:flattened-curvature-reduced}
    Let $M$ be a connected piece of the band of some nested prismoid $P$, which is bounded at the top and base by polygonal paths $\Ga$ and $\Omega$, respectively. 
    Then $\bar{\Ga}$ and $\bar{\Omega}$ are stretchings of $\tilde{\Ga}$ and $\tilde{\Omega}$, respectively.
    The curvature at every interior vertex of $\bar{\Ga}$ and $\bar{\Omega}$ must remain larger than $0$.
\end{observation}

First, recall some definitions from the literature.
The total angle of a vertex of a polytope is the sum of the
incident face angles. 
For any polytope, the total angle is $\leq2\pi$
with equality if and only if the incident faces lie in a plane (which does not
occur in our situation).
Let $M$ be a connected piece of the band of some nested prismatoid $P$ which is bounded by the polygonal path $[a_1,\dots,a_k,b_K,\dots,b_1,a_1]$.
We call the boundary subpaths $[a_1,\dots,a_k]$ and $[b_1,\dots,b_K]$ the top and base boundary of $M$, respectively.

For $i\in\{2,\dots,k-1\}$ (the indices of the interior vertices of the top boundary) let $\alpha_A(a_i)$ be the incident top angle at $a_i$, and  $\alpha_M(a_i)$ the sum of the incident band angles. Similarly, for $j\in\{2,\dots,K-1\}$ (the indices of the interior vertices of the base boundary), let $\beta_B(b_j)$ be the incident base angle at $b_j$, and  $\beta_M(b_j)$ the sum of the incident band angles.
Further, define $\alpha_{\bar{M}}(a_i)$ and $\beta_{\bar{M}}(b_j)$ for the unfolded band piece $\bar{M}$ analogously to $\alpha_M(a_i)$ and $\beta_M(b_j)$.

Now let $\bar{M}$ be an unfolding of $M$ into the plane.
Due to the assumption that $P$ is nested, we have
\[
\beta_B(b_j)\ <\ \beta_M(b_j)\ =\ \beta_{\bar{M}}(b_j)\ <\ \pi\,.
\] 
Since the total angle at $a_i$ is smaller than $2\pi$, we get 
\[
\alpha_A(a_i)\ <\ 2\pi -\alpha_{\bar{M}}(a_i)\ =\ 2\pi -\alpha_M(a_i)\ <\ \pi\ < \alpha_M(a_i)\,.
\]
That is, the unfolding of $[b_1,\dots,b_K]$, the base boundary of $M$, is a stretching of its orthogonal projection to the plane.
Likewise, $[a_1,\dots, a_k]$ is stretched by its unfolding.
This establishes the first part of the observation. 

The curvature at every interior vertex of $\bar{\Ga}$ and $\bar{\Omega}$ must remain larger than $0$ since $P$ is properly nested in the sense that the boundary of $\tilde{A}$ and $B$ do not intersect.
Note that without proper nesting this would not necessarily be the case; see, for example, the band unfolding of a prism in Figure \ref{fig:patterns} (right).

\section{Cutting strategy and placing the top}\label{sec:cutting}
In Section \ref{sec:observations} we collect some observations about band unfoldings and then derive our cutting strategy from these insights.
In particular, we make observations about band pieces whose top and base end edges are parallel.
This is not the case in general. However, in Section \ref{sec:cutting-and-placing} we will
show that there exists a cutting strategy for arbitrary nested prismatoids so that the last facets of the band pieces can (if they are triangles) be embedded into trapezoids in a way that the resulting strictly larger polyhedral surface is unfoldable without self-overlap, implying (overlap-free) unfoldability of the actual (unextended) nested prismatoid.

\subsection{Observations about band unfoldings}\label{sec:observations}
Let $M$ be a connected piece of the band of a nested prismatoid $P$, bounded at the top by $\Ga\define[a_1,\dots,a_k]$ and at the base by $\Om\define[b_1,\dots,b_K]$.
Assume that the pairs of edges $(a_1,a_2)$, $(b_1,b_2)$, and $(a_{k-1},a_k)$, $(b_{K-1},b_K)$ are each contained in a lateral trapezoid and are thus parallel. 
Then by elementary geometry $\Ga$ and $\Om$ both have the same curvature.
Assume that this curvature is $\leq\pi$.
Then, due to Observation \ref{obs:flattened-curvature-reduced}, the curvature of $\bG$ and $\bar{\Om}$ is smaller than $\pi$, while the curvature at every interior vertex of $\bG$ and $\bar{\Om}$ larger than $0$.
Together with the fact that their end edges must be parallel, this implies that $\bar{M}$ does not self-intersect and any line through an edge of $\bar{\Om}$ induces a closed half plane that contains $\bar{M}$.    

Moreover, since the curvature of $\bG$ is smaller than $\pi$, by elementary arithmetic there must exist an edge $\bar{e}\define(\bar{a}_i, \bar{a}_{i+1})$ so that if its relative interior is removed from $\bG$, each of the two remaining subpaths of $\bG$ either consists of a single vertex or has a curvature $\leq\pi/2$ and is thus radially monotone.
Hence, Observations \ref{obs:monotone-intersections} and \ref{obs:flattened-curvature-reduced} imply that if we attach $A$ (thereafter $\bar{A}$) to $\bar{M}$ along $\bar{e}$, it does not intersect $\bar{M}$ anywhere else.
Moreover, any line through an edge of $\bar{\Ga}$ induces a closed half plane that contains $\bar{A}$. In particular, $\bar{A}$ cannot properly intersect the affine rays emanating from $\bar{a}_2$ and $\bar{a}_{k-1}$ through the edges $(\bar{a}_1,\bar{a}_2)$ and $(\bar{a}_{k-1},\bar{a}_k)$, respectively. 
Composing the latter rays with the subpath $[\bar{a}_2,\dots,\bar{a}_{k-1}]$ yields an unbounded curve, say $\bG'$, that intersects $\bar{A}$ only in the edge $(\bar{a}_i, \bar{a}_{i+1})$. 
Since $\bG$ and $\bar{\Om}$ are parallel in their ends, every line through an edge of $\bar{\Om}$ induces a closed half plane whose interior contains $\bG'$ and thus also $\bar{A}$, cf. Figure \ref{fig:barA-and-barM-containment}.
We summarize the relevant aspects of these findings.

\begin{figure}
	\begin{center}
	
	
	\begin{tikzpicture}[scale=0.7]

    
	
	\draw (6,-1) -- (7,0.5) -- (4,0) -- (2,-3) -- (6,-6.5) -- (6,-5) -- (4.5,-4) -- (4,-2) -- cycle;
	
	\draw (6,-1) -- (7,-1) -- (8,0.5) -- (7,0.5);
	\draw (6,-6.5) -- (6.5,-6.5) -- (6.5,-5) -- (6,-5);
	
	\draw (4,-2) -- (4.5,-4) -- (6,-5) -- (6.5,-5) -- (7,-1) -- (6,-1) -- (4,0) -- (4,-2) -- (2,-3) -- (4.5,-4) -- (6,-6.5);
	
	\draw (6.5,-5) -- (7,-1) -- (6,-1) -- (4,-2) -- (4.5,-4) -- (6,-5) -- cycle;

	\draw (3.5,-2.9) node{$\tilde{M}$};
	\draw (4.6,-2.8) node{$\tilde{e}$};
	\draw (5.9,-2.8) node{$\tilde{A}$};
	\draw[ultra thick] (4,-2) -- (4.5,-4);
	\draw (7.45,-1.1) node{$\tilde{a}_1$};
	\draw (6.95,-4.95) node{$\tilde{a}_k$};
	\draw (8.4,0.5) node{$\tilde{b}_1$};
	\draw (7,-6.4) node{$\tilde{b}_K$};

	\end{tikzpicture}
	\hspace{1cm}
		\begin{tikzpicture}[scale=0.7]

    \draw[dotted] (2,0) -- (10,0);
    \draw (4,0) -- (6,-1.41) -- (8,0) -- cycle;
    \draw (4,0) -- (3.92,-2.23) -- (6,-1.41);
    \draw[dotted] (0,-4.11) -- (5.5, 1.54);
    \draw (4,0) -- (1.50,-2.59) -- (3.92,-2.23);
    \draw (1.50,-2.59) -- (3.81,-4.29) -- (3.92,-2.23);
    \draw[dotted] (0.576,-0.678) -- (4.734,-9.282);
    \draw (1.50,-2.59) -- (3.81,-7.37) -- (3.81,-4.29);
    \draw (3.81,-7.37) -- (4.75,-5.83) -- (3.81,-4.29);
    
    \draw[dotted] (4.75+0.39,-5.83-0.3) -- (4.75+1.95,-5.83-1.5);
    \draw[dotted] (3.81-1.95,-7.37+1.5) -- (3.81+1.95,-7.37-1.5);
    \draw (4.75,-5.83) -- (4.75+0.39,-5.83-0.3) -- (3.81+0.39,-7.37-0.3) -- (3.81,-7.37);
    
    \draw[dotted] (5.99,+0.2) -- (10.99,-0.3);
    \draw[dotted] (6.99,-1.51) -- (8.99,-1.71);
    \draw (8,0) -- (8.99,-0.1) -- (6.99, -1.51) -- (6,-1.41);
    
	\draw (5.51154577279, -5.74077289162) -- (6.95693222487, -1.977681342905) -- (5.98615933743, -1.737681340815) -- (3.92,-2.23) -- (3.81,-4.29) -- (5.02615932907, -5.620772890575) -- cycle;
    	
    \draw[ultra thick] (3.92,-2.23) -- (3.81,-4.29);
    \draw (4.21,-3.36) node{$\bar{e}$};
    	
    \draw (3.25,-3) node{$\bar{M}$};
	\draw (5.5,-3.3) node{$\bar{A}$};
	
	\draw (7.5,-1.7) node{$\bar{a}_1$};
	\draw (5.55,-6.2) node{$\bar{a}_k$};
	\draw (9.4,-0.1) node{$\bar{b}_1$};
	\draw (4.75,-7.8) node{$\bar{b}_K$};

	\end{tikzpicture}
	\end{center}
\caption{Projection $\tilde{M}$ and unfolding $\bar{M}$ of band piece $M$.
First and last edges of $\bar{\Omega}\define[\bar{b}_1,\dots,\bar{b}_K]$ and $\bar{\Ga}\define[\bar{a}_1,\dots,\bar{a}_K]$ are parallel.
Dotted lines through every edge of $\bar{\Omega}$ induce halfspaces that contain $\bar{N}$, the surface consisting of $\bar{M}$ and $\bar{A}$.}\label{fig:barA-and-barM-containment}
\end{figure}

\begin{observation}\label{obs:barA-and-barM-containment}
    In the above constellation, where the first and last edges of $\bG$ and $\bar{\Om}$ are parallel, $\bar{M}$ does not self-intersect.
    Moreover, there exists an edge $e$ of $\bG$ so that if we attach $A$ (thereafter $\bar{A}$) to $\bar{M}$ along $e$, the resulting flat polyhedral surface, say $\bar{N}$, does not self-intersect and every line through an edge of $\bar{\Om}$ induces a closed half plane that contains $\bar{N}$.
\end{observation}

\subsection{Cutting and placing the top}\label{sec:cutting-and-placing}
We will now devise a cutting-scheme that recreates the above ideal constellation sufficiently well to harness all its advantages.
Let $P$ be a nested prismatoid with top $A$, base $B$, and corresponding flat prismatoid $P^0$.
We assume that $A$ is an $m$-gon with boundary $[a_1,\dots,a_m]$ and $B$ an $n$-gon with boundary $[b_1,\dots,b_n]$.

Now pick any vertex of $B$. 
By rotating the indices, if necessary, we can assume that we picked $b_1$.  
Let $L_1$ be the line through the edge $(b_n,b_1)$.
We say a line $L$ supports a convex polygon $C$ in a plane $H$ if it lies in $H$, has a nonempty intersection with $C$ and $C$ is contained in one of the two closed half planes induced by $L$. 
Let $L_2$ and $L_3$ be the unique disjoint supporting lines of $\tilde{A}$ which are parallel to $L_1$, and let $L_2$ be the one closer to $L_1$. Further, let $L_4$ be the unique supporting line of $B$ which is parallel to $L_1$ and disjoint from it, cf. Figure \ref{fig:cut-scheme} (left).

We denote by $K$ the smallest index so that $b_K$ is contained in $L_4$.
Further, enumerate the indices of $\tilde{A}$ so that $\tilde{a}_1$ is contained in $L_2$, but $\tilde{a}_2$ is not and denote by $k$ the smallest index so that $\tilde{a}_k$ is contained in $L_3$. Again,  cf. Figure \ref{fig:cut-scheme} (left).
\begin{observation}\label{obs:edges-adjacent-oneinfacet}
    The vertices $a_1$ and $b_1$, as well as $a_k$ and $b_K$ are adjacent.
    Moreover, the edge $(a_1, b_1)$ lies in a lateral facet which contains the edge $(b_n, b_1)$. 
\end{observation}

\begin{figure}
	\begin{center}
	\begin{tikzpicture}[scale=0.75]
	
	\draw[thick,fill=lightgray!20] (6,1) -- (8,2) -- (7.5,4) -- (6,5) -- (4.5,4) -- (4,2) -- cycle;
		
	\draw[fill=lightgray,thick,opacity=1] (8,0) -- (10,3) -- (6,6.5) -- (6,5) -- (7.5,4) -- (8,2) -- (6,1) -- cycle; 
	\draw[fill=lightgray!50,thick,opacity=1] (6,1) -- (8,0) -- (4,0) -- (2,3) -- (6,6.5) -- (6,5) -- (4.5,4) -- (4,2) -- cycle;
	\draw (6,1) -- (4,0) -- (4,2) -- (2,3) -- (4.5,4) -- (6,6.5) -- (7.5,4) -- (10,3) -- (8,2) -- (8,0);
	\draw[ultra thick] (6,6.5) -- (6,5);
	\draw[ultra thick] (8,0) -- (6,1);
	
	\draw (8.2,-.44) node{$b_1$};
	\node at (8,0) [circle,fill,inner sep=1.5pt]{};
	\draw (4,-.44) node{$b_n$};
	\node at (4,0) [circle,fill,inner sep=1.5pt]{};
	\draw (6,1.41) node{$\tilde{a}_1$};
	\node at (6,1) [circle,fill,inner sep=1.5pt]{};
	\draw (6.05,6.95) node{$b_K$};
	\node at (6,6.5) [circle,fill,inner sep=1.5pt]{};
	\draw (6.05,4.55) node{$\tilde{a}_k$};
	\node at (6,5) [circle,fill,inner sep=1.5pt]{};

	\draw (8.6,3) node{$\tilde{M}^+$};
	\draw (3.6,3) node{$\tilde{M}^-$};
	\draw (5.95,3) node{$\tilde{A}$};
	
	\draw[thick,dotted] (1,1) -- (11,1);
	\draw[thick,dotted] (1,0) -- (4,0);
	\draw[thick,dotted] (8,0) -- (11,0);
	\draw[thick,dotted] (1,5) -- (11,5);
	\draw[thick,dotted] (1,6.5) -- (11,6.5);
	
	\draw (1.5,-.44) node{$L_1$};
	\draw (1.5,1.38) node{$L_2$};
	\draw (1.5,6.95) node{$L_4$};
	\draw (1.5,4.55) node{$L_3$};

	\end{tikzpicture}
	\begin{tikzpicture}[scale=0.75]
    

	\draw[fill=lightgray,thick,dotted] (1,0) -- (8,0) -- (6,1) -- (1,1) -- cycle;
	\draw[fill=lightgray,thick,dotted] (1,5) -- (6,5) -- (6,6.5) -- (1,6.5) -- cycle;

	\draw[fill=lightgray,thick,opacity=1] (8,0) -- (10,3) -- (6,6.5) -- (6,5) -- (7.5,4) -- (8,2) -- (6,1) -- cycle; 
	\draw[thick,opacity=1] (6,1) -- (8,0) -- (4,0) -- (2,3) -- (6,6.5) -- (6,5) -- (4.5,4) -- (4,2) -- cycle;
	\draw (6,1) -- (4,0) -- (4,2) -- (2,3) -- (4.5,4) -- (6,6.5) -- (7.5,4) -- (10,3) -- (8,2) -- (8,0);
	\draw[ultra thick] (6,6.5) -- (6,5);
	\draw[ultra thick] (8,0) -- (6,1);
	
	\draw (8.2,-.44) node{$b_1$};
	\node at (8,0) [circle,fill,inner sep=1.5pt]{};
		\draw (4,-.44) node{$b_n$};
	\node at (4,0) [circle,fill,inner sep=1.5pt]{};
	\draw (6,1.41) node{$\tilde{a}_1$};
	\node at (6,1) [circle,fill,inner sep=1.5pt]{};
	\draw (6.05,6.95) node{$b_K$};
	\node at (6,6.5) [circle,fill,inner sep=1.5pt]{};
	\draw (6.05,4.55) node{$\tilde{a}_k$};
	\node at (6,5) [circle,fill,inner sep=1.5pt]{};
	
	\draw (1.5,-.44) node{$g_1$};
	\draw (1.5,1.38) node{$g_2$};
	\draw (1.5,6.95) node{$g_4$};
	\draw (1.5,4.55) node{$g_3$};
	
	
	\draw (8.4,3) node{$\tilde{N}$};
	\draw (1.95,0.5) node{$\tilde{T}_1$};
	\draw (1.95,5.75) node{$\tilde{T}_2$};
	
	
	\end{tikzpicture}
	\end{center}
\caption{Left: Determining the lateral cut edges; the cut edges $(b_K,a_k)$ and $(a_1,b_1)$ are printed fat. Right: Gluing trapezoids to one of the band pieces.}\label{fig:cut-scheme}
\end{figure}

How do we see this? Pick an arbitrary vertex of the base, say $b_i$, and let $H_-$ and $H_+$ be two planes which contain the lateral facets with base edges $(b_{i-1}, b_i)$ and $(b_i,b_{i+1})$, respectively. 
Then the intersections of $H_-$ and $H_+$ with $H_A$ project orthogonally into lines $L_-$ and $L_+$ in $H_B$ which are parallel to $(b_{i-1}, b_i)$ and $(b_i,b_{i+1})$, respectively. 
Relabel $A$ so that $\tilde{a}_1$ is contained in $L_-$, but $\tilde{a}_2$ is not, and let $k$ be the smallest index so that $\tilde{a}_k$ is contained in $L_+$.
Then by construction every lateral edge incident to $b_i$ must contain a vertex of the path $[a_1, \dots, a_k]$, and the first and the last vertices of this path are contained in the lateral facets with base edges $(b_{i-1}, b_i)$ and $(b_i,b_{i+1})$, respectively.
This establishes the claim.

Now cut the lateral edges $(a_1,b_1)$ and $(a_k, b_K)$, as well as all top and base edges.
This dissects $P$ into four pieces, the top $A$, the base $B$ and two band pieces. 
We denote the band piece in counterclockwise direction from $(a_1,b_1)$ by $M^+$ and the one in clockwise direction by $M^-$, cf. Figure \ref{fig:cut-scheme} (left).

Next, we recreate our ideal constellation outlined above by embedding $M^+$ in a strictly larger polyhedral surface that satisfies all conditions which lead to Observation \ref{obs:barA-and-barM-containment}.
To this end, let $g_1$, $\tilde{g}_2$, $\tilde{g}_3$, and $g_4$ be four parallel line segments of nonzero but otherwise arbitrary length with the following properties, cf. Figure \ref{fig:cut-scheme} (right).
\begin{itemize}
    \item They lie in $L_1$, $L_2$, $L_3$, and $L_4$, respectively. 
    \item They originate in $b_1$, $\tilde{a}_1$, $\tilde{a}_k$, and $b_K$, respectively. 
    \item They all extend into the same direction, the one where neither $\tilde{g}_2$, nor $\tilde{g}_3$ intersect the interior of $\tilde{M}^+$.
\end{itemize}
Let $g_2$ and $g_3$ be the orthogonal projections of $\tilde{g}_2$ and $\tilde{g}_3$ to $H_A$ and glue the two trapezoids, say $T_1$ and $T_2$, which arise as the convex hulls of $g_1$ and $g_2$, resp., $g_3$ and $g_4$, to $M^+$ along their common edges $(a_1, b_1)$ and $(a_k, b_K)$, respectively. Again, cf. Figure \ref{fig:cut-scheme} (right). 
This constructed polyhedral surface has a top and base curvature of exactly $\pi$ and both its end facets are trapezoids.
Hence, Observation \ref{obs:barA-and-barM-containment} applies to it -- and thus to its subset $M^+$. 
Applying the analogous construction to $M^-$, we get:
\begin{observation}\label{obs:Mplus-Mminus-contained}
$\bar{M}^+$ does not self-intersect. Moreover, there exists an edge $e$ of its top boundary so that if we attach $A$ (thereafter $\bar{A}$) to $\bar{M}^+$ along $e$, the resulting flat polyhedral surface does not self-intersect and every line through an edge of the base boundary of $\bar{M}^+$ induces a closed half plane that contains it.
Analogous statements hold for $\bar{M}^-$.
\end{observation}

\section{Gluing the band pieces to the base}\label{sec:gluing}
The cutting strategy presented above dissects a nested prismatoid into four pieces, $A$, $B$, $M^+$, and $M^-$. 
Moving on, we attach the top $A$ to $M^+$ along the edge $e$ from Observation \ref{obs:Mplus-Mminus-contained}. 
We will denote the polyhedral surface consisting of $A$ glued to $M^+$ along $e$ by $N$.
There are technical reasons we will explain below to choose $M^+$ and not $M^-$ for the attachment of the top.

We will now single out two edges along which the three pieces $B$, $M^-$ and $N$ can be reattached so that the resulting polyhedral surface whose unfolding does not self-intersect.
We will denote by $e_-$ the edge that connects $M^-$ to $B$ and by $e_+$ the edge  that connects $B$ to $M^+$.
To this end, we distinguish two cases.
\begin{itemize}
    \item \textbf{Case 1:} There exists a vertex of the base $B$ with a curvature $\geq\frac\pi2$.
    \item \textbf{Case 2:} No such vertex exists.
\end{itemize}
We will treat the first case in detail and then show how to reduce the second case to the first.

\begin{figure}
	\begin{center}
	\begin{tikzpicture}[scale=0.6]
	
	\draw[thick,fill=lightgray!20] (6,1) -- (6.7,2) -- (7,5) -- (4.5,4) -- (4,2) -- cycle;
	\draw[fill=lightgray,thick,opacity=1] (8,0) -- (8,6.5) -- (7,5) -- (6.7,2) -- (6,1) -- cycle;
	\draw[thick] (6.7,2) -- (8,0) -- (7,5);
	\draw[fill=lightgray!50,thick,opacity=1] (6,1) -- (8,0) -- (4,0) -- (2.5,3.5) -- (8,6.5) -- (7,5) -- (4.5,4) -- (4,2) -- cycle;
	\draw (6,1) -- (4,0) -- (4,2) -- (2.5,3.5) -- (4.5,4) -- (8,6.5);
	\draw[ultra thick] (8,6.5) -- (7,5);
	\draw[ultra thick] (8,0) -- (6,1);
	
	\draw (8.2,-.44) node{$b_1$};
	\node at (8,0) [circle,fill,inner sep=1.5pt]{};
	\draw (4,-.44) node{$b_n$};
	\node at (4,0) [circle,fill,inner sep=1.5pt]{};
	\draw (5.85,1.63) node{$\tilde{a}_1$};
	\node at (6,1) [circle,fill,inner sep=1.5pt]{};
	\draw (6.68,1.45) node{$\tilde{v}$};
	\node at (6.8,1) [circle,fill,inner sep=1.5pt]{};
	\draw (8.05,6.95) node{$b_K$};
	\node at (8,6.5) [circle,fill,inner sep=1.5pt]{};
	\draw (7.55,4.6) node{$\tilde{a}_k$};
	\node at (7,5) [circle,fill,inner sep=1.5pt]{};

	\draw (8.8,3) node{$\tilde{M}^+$};
	\draw (3.75,3.3) node{$\tilde{M}^-$};
	\draw (5.55,3) node{$\tilde{A}$};
	
	\draw[thick,dotted] (1,1) -- (10,1);
	\draw[thick,dotted] (1,0) -- (4,0);
	\draw[thick,dotted] (8,0) -- (10,0);
	\draw[thick,dotted] (1,5) -- (10,5);
	\draw[thick,dotted] (1,6.5) -- (10,6.5);
	
	\draw (1.5,-.44) node{$L_1$};
	\draw (1.5,1.38) node{$L_2$};
	\draw (1.5,6.95) node{$L_4$};
	\draw (1.5,4.55) node{$L_3$};

	\end{tikzpicture}
	
	\begin{tikzpicture}[scale=0.6]
	
	\draw[thick,dotted] (8,0) -- (8,-7);
	\draw[thick,fill=lightgray!20] (10,1) -- (9.3,2) -- (9,5) -- (11.5,4) -- (12,2) -- cycle;
	\draw[fill=lightgray,thick,opacity=1] (8,0) -- (8,6.5) -- (9,5) -- (9.3,2) -- (10,1) -- cycle;
	\draw[thick] (9.3,2) -- (8,0) -- (9,5);
	\draw[fill=lightgray!50,thick,opacity=1] (6,-1) -- (8,0) -- (4,0) -- (2.5,-3.5) -- (8,-6.5) -- (7,-5) -- (4.5,-4) -- (4,-2) -- cycle;
	\draw (8,0) -- (4,0) -- (2.5,3.5) -- (8,6.5);
	\draw (6,-1) -- (4,0) -- (4,-2) -- (2.5,-3.5) -- (4.5,-4) -- (8,-6.5);
	
	\draw (8.35,-.44) node{$b_1$};
	\draw (8.35,-3.44) node{$g$};
	\draw (4.3,.5) node{$b_n$};
	\draw (6.2,.34) node{$e_-$};
	\draw (7.5,3) node{$e_+$};
	
	\draw[dotted] (6,-1) -- (6.8,-1) -- (8,0);
	\draw (6.2,-1.4) node{$\tilde{a}'_1$};
	\draw (7.15,-1.35) node{$\tilde{v}'$};
	\node at (6.8,-1) [circle,fill,inner sep=1.5pt]{};
	
	\draw (8.05,6.95) node{$b_K$};
	\draw (8.45,-6.5) node{$b'_K$};
	\draw (7.4,-4.75) node{$\tilde{a}'_k$};

	\draw (9.5,6) node{$(\tilde{M}^+)'$};
	\draw (5.25,-3) node{$(\tilde{M}^-)'$};
	\draw (5.7,2.5) node{$B$};
	\draw (10.35,3) node{$\tilde{A}'$};
	\end{tikzpicture}
		\begin{tikzpicture}[scale=0.6]
	
	\draw[thick,dotted] (8,0) -- (8,-9);
	\draw (8,0) -- (8,6.5);



\draw[fill=lightgray!50,thick,opacity=1] (8,0) -- (6.0, -2.23606797749979) -- (4,0) -- cycle;
\draw[fill=lightgray!50,thick,opacity=1] (4,0) -- (3.8426213483334735, -2.8240453183331935) -- (6.0, -2.23606797749979) -- cycle;
\draw[fill=lightgray!50,thick,opacity=1] (4,0) -- (0.9728177830027867, -2.310014680710026) -- (3.8426213483334735, -2.8240453183331935) -- cycle;
\draw[fill=lightgray!50,thick,opacity=1]  (0.9728177830027867, -2.310014680710026) -- (2.7459221581020317, -4.569682991547684) -- (3.8426213483334735, -2.8240453183331935) -- cycle;
\draw[fill=lightgray!50,thick,opacity=1]  (0.9728177830027867, -2.310014680710026) -- (0.11382590869309217, -8.51582912845739) -- (2.7459221581020317, -4.569682991547684) -- cycle;
\draw[fill=lightgray!50,thick,opacity=1]  (0.11382590869309217, -8.51582912845739) -- (2.490368942338184, -7.250110672704696) --  (2.7459221581020317, -4.569682991547684) -- cycle;

\draw[fill=lightgray,thick,opacity=1] (8,6.5) -- (10.23606797749979,5) -- (8,0) -- cycle;
\draw[fill=lightgray,thick,opacity=1] (10.23606797749979,5) -- (10.394388288923338, 1.9891970042876015) -- (8,0) -- cycle;
\draw[fill=lightgray,thick,opacity=1]  (10.394388288923338, 1.9891970042876015) -- (10.872293671420163, 0.8659844485437826) -- (8,0) -- cycle;


	\draw[thick,fill=lightgray!20] (11.04649539474863, 0.957326701133101) -- (10.394388288923338, 1.9891970042876015) -- (10.23606797749979,5) -- (12.68617608937741, 3.8833217830943947) -- (13.091389798001831, 1.8619851336609459) -- cycle;

	\draw (8,0) -- (4,0) -- (2.5,3.5) -- (8,6.5);

	
 	\draw (8.35,-.44) node{$b_1$};
 	\draw (8.35,-3.44) node{$g$};
 	\draw (4.3,.5) node{$b_n$};
 	\draw (6.2,.34) node{$e_-$};
 	\draw (7.5,3) node{$e_+$};
 	\draw (2.9,-7.55) node{$\bar{a}_k$};
 	\draw (0,-8.9) node{$\bar{b}_K$};
 	\draw (8.05,6.95) node{$b_K$};
	\draw[dotted] (6,-2.23606797749979) -- (6.8,-2.23606797749979) -- (8,0);
	\draw (7.05,-2.7) node{$\bar{v}$};
	\node at (6.8,-2.23606797749979) [circle,fill,inner sep=1.5pt]{};
 	\draw (6.2,-2.8) node{$\bar{a}_1$};

 	\draw (9.1,4.7) node{$\bar{M}^+$};
 	\draw (3.15,-1.7) node{$\bar{M}^-$};
 	\draw (5.7,2.5) node{$B$};
 	\draw (11.45,2.9) node{$\bar{A}$};
	\end{tikzpicture}
	\end{center}
\caption{\textbf{Case 1:} Projecting and cutting (top), reflecting (bottom left) and unfolding (bottom right).}\label{fig:case1-unfolding}
\end{figure}

\textbf{Case 1:} Label $B$ so that the vertex with curvature $\geq\frac\pi2$ is $b_1$ and apply our cutting strategy, cf. Figure \ref{fig:case1-unfolding} (top). 
Then set $e_-\define(b_n,b_1)$ and $e_+\define(b_1,b_2)$, cf. Figure \ref{fig:case1-unfolding} (bottom left and right).
By Observation \ref{obs:Mplus-Mminus-contained} we can establish unfoldability in the first case by proving that $\bar{M}^-$ intersects the line through $e_+$ nowhere except in $b_1$.
Let $g$ be the outward normal ray of $e_-$ at $b_1$.
Since the curvature at $b_1$ is $\geq\pi/2$, we are done if we can prove that $\bar{M}^-$ intersects $g$ in $b_1$ and nowhere else.
 For this it suffices to prove that the base boundary of $\bar{M}^-$ intersects $g$ only in $b_1$ and its top boundary does not intersect $g$ at all.

By Observation \ref{obs:edges-adjacent-oneinfacet} our cutting strategy ensures the following: If there is more than one lateral edge incident to $b_1$, then $(a_1,b_1)$ is the first one of them counted in counterclockwise direction and there must exist a lateral facet that contains the edges $e_-=(b_n, b_1)$ \emph{and} $(a_1,b_1)$. (This is not the case for $M^+$. Hence, our choice above.)
This facet is either a trapezoid or a triangle. 
If it is not a trapezoid, embed the facet in a trapezoid by taking its convex hull with a point, say $v$, that lies on $L_2$ and projects into the interior of $B$, but not into $\tilde{M}^-$, cf. Figure \ref{fig:case1-unfolding} (all subfigures).
By construction, this does not increase the top curvature of $\tilde{M}^-$ beyond $\pi$.
We will thus assume without loss of generality that the lateral facet containing $e_-$ is a trapezoid.

\begin{figure}
	\begin{center}	
	\begin{tikzpicture}[baseline=(current bounding box.north),scale=0.7]
	
	\begin{scope}
	\clip (-5,0) rectangle (5,5);
	\draw[dotted] (0,0) circle(5);
	\end{scope}
	\begin{scope}
	\clip (-3,1) rectangle (7,6);
	\draw[dotted] (2,1) circle(5);
	\end{scope}
	
	\draw[thick, fill=gray!20] (-3,1) -- (-5,0) -- (-2,-1) -- (4,2) -- cycle;
	
	\draw[dotted, thick] (5.53553390593,1) -- (3.53553390593,0) -- (-2,-1) -- (4,2) -- cycle;
	\draw[fill=gray!10] (-2,-1) -- (-3.53553390593,0) -- (-1.53553390593,1) -- (4,2)  -- cycle;
	\draw[thick, dotted] (3.53553390593,3.53553390593) -- (3.53553390593,0);
	\draw[thick, dotted] (5.53553390593,4.53553390593) -- (5.53553390593,1);

	\draw[thick, dotted] (-8,0) -- (3.53553390593,0); 
	\draw[thick, dotted] (-6,1) -- (5.53553390593,1); 
	\draw[thick, fill=gray!30,opacity=0.5] (-2,-1) -- (3.53553390593,3.53553390593) -- (5.53553390593,4.53553390593) -- (4,2) -- cycle;

	\draw[thick,dotted] (4,2) -- (-2,-1);
	\draw[thick,dotted] (-9,-1) -- (-2,-1);
	\draw[thick] (7,-0.75) -- (-2,-1);
	\draw (2.8,-1.3) node{$e_+$};
	\draw (-8.5,-1.4) node{$g$};
    \draw (3.75,-.45) node{$\tilde{a}_1$};
    \draw (3.6,3.95) node{$a_1$};
    \draw (4.95,2.2) node{$b_n$};
    \draw (1.8,0.4) node{$e_-$};
    \draw (-1.8,-1.4) node{$b_1$};
    \draw (-4.95,-.5) node{$\bar{a}_1$};
	\end{tikzpicture}
	\end{center}
	\caption{Projection, reflection of projection, and unfolding of lateral trapezoid incident to $e_-\define(b_n, b_1)$ in \textbf{Case 1}.  }\label{fig:unfolding-vs-projection}
\end{figure}

 Since $P$ is properly nested, i.e., $\tilde{A}$ is contained in the interior of $B$, and due to the curvature at $b_1$, the reflection of $\tilde{M}^-$ at $e_-$ can intersect $g$ in at most two points. These are the endpoints of its base boundary, $b_1$ and the reflection of $\tilde{b}_\ell$ at $e_-$. 
 Let $p$ be a point in either the top or the base boundary of $M^-$, and $p'$ the reflection of $\tilde{p}$ at $e_-$. 
 We will show that the distance of $\bar{p}$ to $g$ is larger than or equal to the distance of $p'$ to $g$, with equality if and only if $p$ is contained in $e_-=(b_n,b_1)$ or $(a_m,a_1)$, that is, if it lies in the top or base edge of the lateral trapezoid that contains $e_-$.

 If $p$ lies in $e_-$, this is clear, because $e_-=\tilde{e}_-=\bar{e}_-$.
 By elementary geometry, the edges $(\tilde{a}_n,\tilde{a}_1)$, its reflection at $e_-$, and $(\bar{a}_n,\bar{a}_1)$ are parallel.
 Moreover, $\tilde{a}_1$, its reflection at $e_-$, and $\bar{a}_1$ lie on a line which is perpendicular to $e_-$, cf. Figure \ref{fig:unfolding-vs-projection}.
 Hence, $(\tilde{a}_n,\tilde{a}_1)$, its reflection at $e_-$, and $(\bar{a}_n,\bar{a}_1)$ can be translated into each other by translating them parallel to $g$.
 This proves the claim for $p\in(a_n,a_1)$.

\begin{figure}
	\begin{center}
	\begin{tikzpicture}[scale=0.5]

    \draw (4,0) -- (8,0) -- (10,3) -- (6,6.5) -- (2,3) -- cycle;
    
    \draw[dotted] (10,3) -- (14,-0.5);
    \draw[dotted] (8,0) -- (14,0);
    \draw[dotted] (13.42,-7.37) -- (13.42,0);
	\node at (13.42,0) [circle,fill,inner sep=1.5pt]{};
	
	\draw[fill=lightgray!50,opacity=1] (6,-1) -- (8,0) -- (4,0) -- (2,-3) -- (6,-6.5) -- (6,-5) -- (4.5,-4) -- (4,-2) -- cycle;
	
	\draw (6,-1) -- (4,0) -- (4,-2) -- (2,-3) -- (4.5,-4) -- (6,-6.5);
	
	\draw (8.5,-0.2) node{$b_1$};
    \draw (3.5,0.3) node{$b_n$};
    \draw (6,0.3) node{$e_-$};
    \draw (7.8,4.2) node{$e_+$};
	\draw (6.6,7.3) node{$(\tilde{M}^+)'$};
	\draw (5.5,-2.8) node{$(\tilde{M}^-)'$};
	\draw (6,3) node{$B$};
	
	\draw[fill=lightgray] (10,3) -- (12.7,5.4) -- (11.45,7.23) -- (10.73,5.11) -- (8.68,5.35) -- (7.48,6.69) -- (6,6.5) -- cycle; 
	\draw (6,6.5) -- (8.68,5.35) -- (10,3) -- (10.73,5.11) -- (12.7,5.4);
	\draw[fill=lightgray!20] (11.45,7.23) -- (10.19,9.08) -- (8.28,8.32) -- (7.48,6.69) -- (8.68,5.35) -- (10.73,5.11) -- cycle;
	
	\draw[ultra thick] (8.68,5.35) -- (10.73,5.11);
	
	\draw (9.6,7) node{$\tilde{A}'$};
	\draw (13,-7) node{$g$};
	\draw (13.65,0.5) node{$p$};

	\end{tikzpicture}
	\begin{tikzpicture}[scale=0.5]

    \draw (4,0) -- (8,0) -- (10,3) -- (6,6.5) -- (2,3) -- cycle;
    \draw (8.5,-0.2) node{$b_1$};
    \draw (3.5,0.3) node{$b_n$};
    \draw (6,0.3) node{$e_-$};
    \draw (7.8,4.2) node{$e_+$};
	\draw (6.9,7.4) node{$\bar{M}^+$};
	\draw (4.7,-3.3) node{$\bar{M}^-$};
    \draw (6,3) node{$B$};
    
    \draw[fill=lightgray!50] (4,0) -- (6,-1.41) -- (8,0);
    \draw[fill=lightgray!50] (4,0) -- (3.92,-2.23) -- (6,-1.41);
    \draw[fill=lightgray!50] (4,0) -- (1.49,-2.59) -- (3.92,-2.23);
    \draw[fill=lightgray!50] (1.49,-2.59) -- (3.81,-4.29) -- (3.92,-2.23);
    \draw[fill=lightgray!50] (1.49,-2.59) -- (3.81,-7.37) -- (3.81,-4.29);
    \draw[fill=lightgray!50] (3.81,-7.37) -- (4.75,-5.83) -- (3.81,-4.29);
    
    \draw[fill=lightgray] (10,3) -- (8.97,5.68) -- (6,6.5);
    \draw[fill=lightgray] (8.97,5.68) -- (7.73,6.99) -- (6,6.5);
    \draw[fill=lightgray] (10,3) -- (10.98,5.24) -- (8.97,5.68);
    \draw[fill=lightgray] (10,3) -- (13.16,4.74) -- (10.98,5.24);
    \draw[fill=lightgray] (13.16,4.74) -- (12.31,7.04) -- (10.98,5.24);
    \draw[fill=lightgray!20] (10.98,5.24) -- (8.97,5.68) -- (7.92,7.14) -- (8.88,8.68) -- (10.86,9.24) -- (11.93,7.27) -- cycle;
    
    \draw (10,7.3) node{$\bar{A}$};

	\end{tikzpicture}
	\end{center}
\caption{Left: Unfolding of $P^0$ according to \textbf{Case 2}. $\tilde{A}'$, the reflection of $\tilde{A}$, is attached at the fat edge of the reflection of $\tilde{M}^+$. Right: Unfolding of $P$}\label{fig:flat-case2}
\end{figure}

 Now denote by $\Om$ the base boundary of $M^-$.
 We can transform the reflection of $\tilde{\Om}$ at $e_-$, which we denote by $\tilde{\Om}'$, into $\bar{\Om}$ as follows. 
 First, decrease the curvature at $\tilde{b}'_{K+1}$ so that it matches the curvature at $\bar{b}_{K+1}$ in $\bar{\Om}$.
 This rotates the edge $(\tilde{b}'_K,\tilde{b}'_{K+1})$ about $\bar{b}_{K+1}$ in clockwise direction (seen from a vantage point above $H_B$). 
 Since the curvature of $\tilde{\Om}'$ is $\leq\pi$ and $e_-$ is perpendicular to $g$, this rotation must increase the distance to $g$ for all points in $(\tilde{b}'_K,\tilde{b}'_{K+1})$, except $\tilde{b}'_{K+1}$. 
 Successively applying this procedure to $\tilde{b}'_{K+2}$,\dots,$\tilde{b}'_n$ yields the claim. 
 A similar argument can be made for the image of the top boundary with one extra step: First perform the rotations, as above. This yields a curve congruent to $\bG$. Then shift this curve parallel to $g$ into $\bG$, which does not change the distance to $g$ for any point in the translated curve. 
 This completes our proof for the first case that there exists a base vertex with curvature $\geq\pi/2$.
 
\textbf{Case 2:} In this case $B$ has no vertex with a curvature $\geq\pi/2$. 
 Then by elementary arithmetic there must exist an index $i$ so that the path $(b_1, b_i)$ has a curvature in $[\frac\pi2,\pi)$.
 We set $e_+\define(b_{i-1}, b_i)$.
 Further, set $e_-\define(b_n,b_1)$, as in the first case.
 Let $L_-$ be the line through $e_-$, $L_+$ the line through $e_+$, and $p$ their intersection. 
 Then the curvature of the polygonal path $[b_1,p,b_{i-1}]$ at $p$ is $\geq\pi/2$, and we can prove the claim that $\bar{M}^-$ does not intersect the outward normal ray $g$ of $(b_1,p)$ emanating from $p$ in full analogy to the above claim that $\bar{M}^-$ intersects $g$ only in $b_1$, where $b_1$ is a base vertex with curvature $\geq\pi/2$.
  This completes our proof for the second case that there exists no base vertex with curvature $\geq\pi/2$, and thus completes the proof of Theorem \ref{thm:nested}.

\section*{Acknowledgments}
The author wants to thank his advisor Michael Joswig for pointing him to D\"urer's Problem and giving him the time to work things out. 
He further wants to thank Joseph O'Rourke for constructive discussions on the subject.
Finally, he wants to thank the referees of this manuscript for suggestions that helped to significantly improve its clarity.

 \bibliographystyle{alpha}
\bibliography{references}
\end{document}